\documentclass[12pt]{amsart}
\usepackage{verbatim}
\usepackage{amsthm}
\usepackage{enumerate}

\newif\ifdeveloping

\let\QED\qed
\newcommand{\prlabel}[1]{\renewcommand{\qed}{\QED${}_{\ref{#1}}$}}
\newcommand{\prtxtlabel}[1]{\renewcommand{\qed}{\QED${}_{{\mbox{\tiny #1}}}$}}
\newcommand{\prnolabel}{\prtxtlabel{}}
\newtheorem{theorem}{Theorem}[section]
\newtheorem{proposition}[theorem]{Proposition}
\newtheorem{lemma}[theorem]{Lemma}
\newtheorem{corollary}[theorem]{Corollary}
\newtheorem{fact}[theorem]{Fact}

\newtheorem{claim}{Claim}[theorem]

\newtheorem{qtheorem}{Theorem}

\theoremstyle{definition}

\theoremstyle{remark}
\newtheorem{remark}{Remark}  
{}

\ifdeveloping
\usepackage[notref,notcite]{showkeys}
\fi

\newcommand{\prtime}{{\count0=\time\divide\count0 by 60
\count1=-\count0\multiply\count1 by 60
\advance\count1 by \time
\the\count0:\the\count1}
}

\def\myheads#1;#2;{
\pagestyle{myheadings}
\markboth{{\sc\hfill #1\hfill\protect\makebox[0cm][r]{\rm\today; \prtime}}}
{{\sc\protect\makebox[0cm][l]{\rm\today;\ \prtime}\hfill #2\hfill}}
\thispagestyle{myheadings}
}

\newcommand{\acal}{{\mathcal A}}
\newcommand{\bcal}{{\mathcal B}}
\newcommand{\ccal}{{\mathcal C}}
\newcommand{\dcal}{{\mathcal D}}

\newcommand{\fcal}{{\mathcal F}}

\newcommand{\hcal}{{\mathcal H}}

\newcommand{\ncal}{{\mathcal N}}

\newcommand{\scal}{{\mathcal S}}
\newcommand{\tcal}{{\mathcal T}}
\newcommand{\ucal}{{\mathcal U}}

\newcommand{\setm}{\setminus}
\newcommand{\empt}{\emptyset}
\newcommand{\subs}{\subset}

\newcommand{\oo}{{{\omega}_1}}

\def\<{\left\langle}
\def\>{\right\rangle}

\def\oo{\omega_1}
\def\br#1;#2;{\bigl[ {#1} \bigr]^ {#2} }
\def\bc#1;#2;{\bigl( {#1} \bigr)^ {#2} }

\def\to{\longrightarrow}

\theoremstyle{plain}

\begin{document}

\author[L. Soukup]{
Lajos Soukup }
\address{
Alfr{\'e}d R{\'e}nyi Institute of Mathematics }
\email{soukup@renyi.hu}

\title{Nagata's conjecture and countably compactifications in generic extensions %
}
\thanks{The preparation of this paper was supported by the
Hungarian National Foundation for Scientific Research grant no. 61600}

\begin{abstract}
 Nagata  conjectured that
every $M$-space is homeomorphic to a closed subspace of
the product of a countably compact space and a metric space.
This conjecture was refuted by Burke and  van Douwen,
and  A.
Kato, independently. 

However, we can show that  
there is a c.c.c. poset $P$ of size $2^{\omega}$ 
such that in $V^P$
Nagata's conjecture holds  
for  each  first countable regular space  from the ground model (i.e. 
if a first countable regular space $X\in V$ is an $M$-space in $V^P$ then 
it is homeomorphic to a closed subspace of
the product of a countably compact space and a metric space in $V^P$).
In fact, we show 
that
every first countable regular space from the ground model   has a first countable
countably compact extension in $V^P$,  and then apply some results of Morita.
As a corollary, we obtain that   every first countable regular space from the ground model
has a maximal first countable extension in model $V^P$.
\end{abstract}

\maketitle \ifdeveloping
\myheads{Countably-compactification};{Countably-compactification};
\fi

\section{Introduction}

A topological space $X$ is called an {\em M-space} (see \cite{N})
if there is a countable collection
of {open covers} $\{\ucal_n:n\in {\omega}\}$ of X, such that:
\begin{enumerate}[(i)]
\item $\ucal_{n+1}$ {star-refines } $\ucal_n$, for all n.
\item If $x_n\in St(x,\ucal_n)$, for all $n$, then the set
$\{x_n : n \in {\omega}\}$ has an accumulation point.
\end{enumerate}

 Nagata, \cite{N}, conjectured that
every $M$-space is homeomorphic to a closed subspace of
the product of a countably compact space and a metric space.

To attack this problem the notion of  countably-compactifiable  spaces
was introduced and studied in \cite{Mo}. A space $X$ is {\em
countably-compactifiable} if it has a {\em
countably-compactification}, i.e. there exists a countably compact
space $S$ such that (1) $X$ is a dense subspace of $S$, and (2)
every countably compact closed subset of $X$ is closed in $S$.
\begin{qtheorem}[Morita, \cite{Mo}]
 An $M$-space satisfies Nagata's conjecture 
(i.e. it is homeomorphic to a closed subspace of
the product of a countably compact space and a metric space) if and only if
it is countably-compactifiable.  
\end{qtheorem}

Burke and  van Douwen in \cite{BvD}, and independently
Kato in \cite{K}  showed that 
there are normal,
first countable
$M$-spaces  which are not countably-compactifiable, 
hence  Nagata's conjecture was refuted.

A countably compact space  is not necessarily a countably-compacti\-fication
of a dense subspace, but 
since a countably compact subspace of a first countable space is closed  
we have
\begin{fact}
A first countable countably compact space  is 
a countably-compacti\-fication
of a dense subspace.
\end{fact}


A first-countable  space $Y$ is said to be a {\em maximal
first-countable extension} of a space $X$ provided
$X$ is a dense subspace of $Y$ and $Y$ is closed  in
any first countable space $Z\supset Y$.

In \cite{TT} the authors considered which
first-countable spaces have first-countable
maximal extensions and whether all do.
They  gave three first-countable spaces without maximal
first-countable extensions.

Since a countably compact subspace $Y$  of a first countable space  $Z$ is closed
in $Z$ we have that 
\begin{fact}
A first-countable countably-compactification $Y$ of a space $X$  is
a maximal first-countable extension of $X$.
\end{fact}

So if you want to construct  maximal first-countable extensions
or  countably compactifications of first countable spaces
the following plan  seems to be natural:
{\em Embed  the first countable spaces into
first countable, countably compact spaces!}  

Although examples from \cite{K} and \cite{TT} are really
sophisticated it is easy to construct a ZFC example of a first
countable space which can not be embedded into a first countable,
countably compact space:

\begin{proposition}\label{pro:psi}
A $\Psi$-space  does not have a first-countable
countably compact extension.
\end{proposition}

\begin{proof}
The underlying set of a $\Psi$-space $X$ is
${\omega}\cup \{x_A:A\in \acal\}$,
where $\acal$ is a maximal almost disjoint family on ${\omega}$, and 
$A$ converges to $x_A$ in $X$ for $A\in\acal$.

Assume on the contrary that  a first countable, countably
compact space $Y$ contains $X$ as a dense subspace.
Let $\{A_n:n\in {\omega}\}$ be distinct elements of $\acal$.
Then $\{x_{A_n}:n\in {\omega}\}$ has an accumulation point $d$ in $Y$.
Since $\{x_{A_n}:n\in {\omega}\}$ is closed in $X$ we have $d\in Y\setm X$.
Since ${\omega}$ is dense in $X$, and so in $Y$, as well, there is a
sequence $D=\{d_n:n\in {\omega}\}\subs {\omega}$ converging to $d$ in $Y$
because $Y$ is first-countable.
 But $\acal$ was maximal so there is $A\in \acal$ with
$|D\cap A|={\omega}$. Hence $x_A$ is an accumulation point of $D$ in $Y$
and so $d=x_A$ because $Y$ is $T_2$. Contradiction.
\end{proof}

The cardinality of a $\Psi$-space is at least $\mathfrak{a}$.
In theorem \ref{tm:s} we show
 that under  Martin's Axiom
every  first countable regular space
 of cardinality $<\mathfrak{c}$  can be embedded,
as a dense subspace, into   a first countable countably compact
regular space. Hence, under Martin's Axiom,  Nagata's conjecture
holds for first countable regular spaces of size less than $\mathfrak c$.

However, the situation changes dramatically if we want to find
a first countable countably compact extension of $X$ in some generic
of the ground model: in theorem \ref{tm:main} we show 
there is a c.c.c. poset $P$ of size $2^{\omega}$ that
every first countable regular space from the ground model   has a first countable
countably compact extension in $V^P$.
Hence, by Corollary \ref{cor:main}, in $V^P$
Nagata's conjecture holds  
for  each  first countable regular space  from the ground model.

\section{First countable, countably compact extensions}

\begin{theorem}
\label{tm:s}
If $\mathfrak{b}=\mathfrak{s}=\mathfrak{c}$ then
every  first countable regular space
 of cardinality $<\mathfrak{c}$  can be embedded,
as a dense subspace, into   a first countable countably compact
regular space.
\end{theorem}

By Morita's results from \cite{Mo} and \cite{Mo2}
Theorem \ref{tm:s}  yields the following corollary:
\begin{corollary}
If Martin's Axiom holds  then Nagata's
conjecture holds for every first countable regular space
$X$ of cardinality $<\mathfrak c$,
i.e. if $X$ is an $M$-space,
then  $X$ is homeomorphic to a closed subspace
of the product of a countably compact space and a metric space.
 Moreover, $X$ has maximal first-countable
extension.

\end{corollary}

\begin{proof}[Proof of Theorem \ref{tm:s}]\prlabel{tm:s}
The proof is based on the following lemma.

\begin{lemma}\label{lm:s}
 Let $X$ be
a first countable regular space $X$, let $A$ be a countable closed discrete
subset of $X$,
and for each $x\in X$ let $\{U(x,n):n<{\omega}\}$
be a neighbourhood base of $x$ such that
$\overline{U(x,n+1)}\subs U(x,n)$.
If $|X|<\min(\mathfrak{b},\mathfrak{s})$ then
there is a first countable regular space $Y$, $|Y|=|X|$,
and for each $y\in Y$ there is a neighbourhood base
$\{U'(y,n):n<{\omega}\}$
of $y$ with
$\overline{U'(y,n+1)}\subs U(y',n)$
such that
\begin{enumerate}[(i)]
\item \label{ext_s} $U'(x,n)\cap X=U(x,n)$ for $x\in X$ and $n<{\omega}$,
\item \label{disj_s} if $U(x,n)\cap U(y,m)=\empt$
then $U'(x,n)\cap U'(y,m)=\empt$,
\item \label{subs_s} if $U(x,n)\subs U(y,m)$ then $U'(x,n)\subs U'(y,m)$,
\item   \label{acc_s} $A$ has an accumulation point in $Y$.
\end{enumerate}
\end{lemma}

\begin{proof}[Proof of Lemma \ref{lm:s}]\prlabel{lm:main}

Since $|X|<\mathfrak{s}$ the family
$\{U(x,n)\cap A:x\in X,n<{\omega}\}$ is not a splitting family, i.e.
there is $D\in \br A;{\omega};$
 such that for each $x\in X$ and for each $n<{\omega}$
either $D_B\subs^* U(x,n)$ or $D_B\cap U(x,n)$ is finite.

Let  $Y=X\cup \{y\}$, where $y$ is a  new point.

We will define the topology of $Y$
as follows.

Let $\theta$ be a large enough regular cardinal and let
$\ncal$ be an elementary submodel of $\<\hcal(\theta),\in,\le\>$
such that $|\ncal|=|X|$, $X, \{U(x,n):x\in X, n\in {\omega}\}\in X$
and $|X|\subs \ncal$. Since $|\ncal|=|X|<\mathfrak{b}$
there is a function $d\in {}^{\omega}{\omega}$
 dominating $\ncal\cap {}^{\omega}{\omega}$.
For $x\in X$ and $n<{\omega}$ let
\begin{displaymath}
U'(x,n)=\left\{
\begin{array}{ll}
U(x,n)\cup\{y\}&\text{ if $D\subs^* U(x,n)$}\\
U(x,n)&\text{ if $D\cap U(x,n)$ is finite}\}
\end{array}
\right..
\end{displaymath}

\newcommand{\dvec}{\vec{D}}
\newcommand{\evec}{\vec{E}}

Let $\dvec$ be a 1-1 enumeration of $D$, and
for $n\in {\omega}$
let
\begin{displaymath}
U'(y,n)=\{y\}\cup\bigcup\{U'(\dvec(k), d(k)+n):n\le k<{\omega}\}.
\end{displaymath}

The family $\{U'(y,n):y\in Y, n\in {\omega}\}$
clearly satisfies (\ref{ext_s}), (\ref{disj_s}) and (\ref{subs_s}).

We intend to define the topology on $Y$ as the one induced
 by the neighbourhood base
$\{U'(y,n):y\in Y, n\in {\omega}\}$.
Next  we should prove that $\{U'(y,n):y\in Y, n\in {\omega}\}$ is a
neighbourhood base of a topology.

\begin{claim}
If $t\in U'(v,n)$ then there is $m$ such that
$U'(t,m)\subs U'(v,n)$.
\end{claim}

\begin{proof}[Proof of the claim]\prnolabel\mbox{ }\\
{\em Case 1: $v,t\in X$.}\\ Then there is $m$ such that
$U(t,m)\subs U(v,n)$. Then $U'(t,m)\subs U'(v,n)$.\\
{\em Case 2: $t\in X$ and $v=y$.}\\
Then $t\in U'(\dvec(k), d(k)+n)$ for some  $k\ge n$ and so
$t\in U(\dvec(k), d(k)+n)$. Thus there is $m$
such that $U(t,m)\subs U(\dvec(k),d(k)+n)$ and so
$U'(t,m)\subs U'(\dvec(k),d(k)+n)\subs U'(v,n)$.\\
{\em Case 3: $t=y$ and $v \in  X$.}\\
Then $D\subs^* U(v,n)$.  Fix $k$ such that
$\dvec(k')\in U(v,n)$ for $k'\ge k$. Pick a function
$g:{\omega}\setm k\to {\omega}$ from $\ncal$
such that $U(\dvec(k'),g(k'))\subs U(v,n)$.
Then there is $m\ge k$ such that
$d(m')\ge g(m')$ for $m'\ge m$. Then
$U(\dvec(m'), d(m')+m)\subs U(v,n)$ for $m'\ge m$, hence
$U'(\dvec(m'), d(m')+m)\subs U'(v,n)$ for $m'\ge m$,
 and so
$U'(y,m)\subs U'(v,n)$.\\
\end{proof}

Hence the family
$\{U'(y,n):y\in Y, n\in {\omega}\}$ can be considered as the
neighbourhood base of a topology  on $Y$.

\begin{claim}\label{t2reg}
If $t\notin U'(v,n)$ then there is $m$ such that
$U'(t,m)\cap U'(v,n+1)=\empt$.
\end{claim}

\begin{proof}[Proof of the claim]\prnolabel\mbox{ }\\
{\em Case 1: $v,t\in X$.}\\
Since $\overline{U(v,n+1)}\subs U(v,n)$  there is $m$ such that
$U(t,m)\cap U(v,n+1)=\empt$. Then
$U'(t,m)\cap U'(v,n+1)=\empt$.\\
{\em Case 2: $t=y\in Y\setm X$ and $v\in X$.  }\\
Since $y_D\notin U'(v,n)$ we have that $D\cap U'(v,n)$ is finite and so
there is $\ell$ such that
$\{\dvec(i):i\ge \ell\}\cap U(v,n)=\empt$.
Then $\{\dvec(i):i\ge \ell\}\cap\overline{U(v,n+1)}=\empt$ and so
we can find a function $g:{\omega}\setm \ell\to {\omega}$
in $\ncal$
such that $ U(\dvec(i),g(i))\cap U(v,n+1)=\empt$ for $i\ge \ell$.
Then there is $m\ge \ell$ such that
$d(m')\ge g(m')$ for $m'\ge m$. Then
$U(\dvec(i),d(i))\cap U(v,n+1)=\empt$ and so
$U'(\dvec(i),d(i))\cap U'(v,n+1)=\empt$ as well for $i\ge m$.
Thus $U'(t,m)\cap U'(v,n+1) =\empt$.\\
{\em Case 3: $t\in X$ and $v=y\in Y\setm X$.  }\\
Since $D$ is closed discrete and $t\notin U(v,n)$ there is $\ell$ such that
$\overline{U(t,\ell)}\cap \{\dvec(i):i\ge n\}=\empt$.
Fix a function $g:{\omega}\setm n\to {\omega}$ in $\ncal $such that
$U(t,\ell)\cap U(\dvec(i),g(i))=\empt$ for $i\ge n$.
There is $k\ge n$ such that $d(i)\ge g(i)$ for $i\ge k$.

Since $t\notin U'(v,n)$ we have $t\notin U(\dvec(i),d(i)+n)$
for $i\ge n$. Thus $t\notin \overline{U(\dvec(i),d(i)+n+1)}$
for $i\ge n$.
Fix $m\ge \ell$ such that
$U(t,m)\cap \cup\{U(\dvec(i),d(i)+n+1):n\le i<\ell\}=\empt$.

Then $U(t,m)\cap U'(v,n+1)=\empt$ and so $U'(t,m)\cap U'(y,n+1)=\empt$
as well.\end{proof}

Since $\bigcap \{U'(x,n):n<{\omega}\}=\{x\}$ claim \ref{t2reg} implies
that $X$ is a regular space.
\end{proof}

Using  the lemma we can easily  prove the theorem.

Let $X$ be a regular first countable space having cardinality
$<\mathfrak{c}=\min(\mathfrak{b},\mathfrak{s})$.
For each $x\in X$ let $\{U(x,n):n<{\omega}\}$
be a neighbourhood base of $x$ such that
$\overline{U(x,n+1)}\subs U(x,n)$.

For ${\alpha}\le 2^{\omega}$
we will construct first countable spaces  $X_{\alpha}$
with bases $\{U_{\alpha}(x,n):x\in X_{\alpha},n<{\omega}\}$
satisfying $\overline{U_{{\alpha}}(x,n+1)}\subs U_{\alpha}(x,n)$,
and sets $A_{\alpha}\in \br X_{\alpha};{\omega};$
such that
\begin{enumerate}[(1)]
\item $X_0=X$, $U_0(x,n)=U(x,n)$
\item $|X_{\alpha}|\le |X|+|{\alpha}|$.
\item \label{ext} $U_{\beta}(x,n)\cap X_{\alpha}=U_{\alpha}(x,n)$ for
${\alpha}<{\beta}$,
$x\in X_{\alpha}$ and $n<{\omega}$,
\item \label{disj} if $U_{\alpha}(x,n)\cap U_{\alpha}(y,m)=\empt$
then $U_{\beta}(x,n)\cap U_{\beta}(y,m)=\empt$ for ${\alpha}<{\beta}$,
\item \label{subs} if $U_{\alpha}(x,n)\subs U_{\alpha}(y,m)$
then $U_{\beta}(x,n)\subs U_{\beta}(y,m)$,
\item   $A_{\alpha}$
has an accumulation point in $X_{{\alpha}+1}$.
\item \label{book} $\{A_{\alpha}:{\alpha}<2^{\omega}\}=\br X_{2^{\omega}};{\omega};$
\end{enumerate}

The construction is straightforward: if ${\alpha}={\beta}+1$ then apply lemma \ref{lm:s}
for $X_{\beta}$, $\{U_{\beta}(x,n):x\in X_{\beta},n<{\omega}\}$ and $A_{\beta}$,
 and for limit ${\alpha}$ take
$X_{\alpha}=\bigcup \{X_{\zeta}:{\zeta}<{\alpha}\}$ and
$U_{\alpha}(x,n)=\bigcup\{U_{\zeta}(x,n):x\in X_{\zeta}\}$. 

By (\ref{book}) every countable subset of $Y=X_{2^{\omega}}$
appears in some intermediate step so
the space $Y$ will be countably compact.
\end{proof}

\begin{theorem}\label{tm:main}
There is a c.c.c. poset of size $2^{\omega}$ such that every first
countable regular space $X$ from the ground model can be embedded,
as a dense subspace, into   a first countable countably compact
regular space $Y$ from the generic extension, and so $X$ has a
countably-compactification in the generic extension.
\end{theorem}

By Morita's results from \cite{Mo} and \cite{Mo2}
Theorem \ref{tm:main}  above yields immediately the following corollary:

\begin{corollary}\label{cor:main}
 There is a c.c.c. poset $P$ of size $2^{\omega}$ such that for
every first countable regular space $X$ from the ground model $V$
Nagata's conjecture holds for $X$ in $V^P$,
i.e. the following holds in $V^P$: if
 $X$ is an $M$-space,
then $X$ is homeomorphic to a closed subspace of the product
of a countably compact space and a metric space.
Moreover, $X$ has  maximal first-countable extension in $V^P$. 
\end{corollary}

\begin{proof}\prlabel{tm:main}
The proof is based on the following lemma.

\begin{lemma}\label{lm:main}
Let $Q=\ccal*R_{\fcal}*\dcal$, where $\ccal$ is the Cohen-poset,
$\fcal$ is a non-principal ultrafilter on ${\omega}$ in $V^\ccal$,
$R_{\fcal}$ introduces a pseudo intersection of the elements of $\fcal$,
and
$\dcal$ is the standard c.c.c poset which adds a dominating real
to $V^{\ccal*R_U}$. Let $X$ be
a first countable regular space from the ground model $V$,
and for each $x\in X$ let $\{U(x,n):n<{\omega}\}$
be a neighbourhood base of $x$ such that
$\overline{U(x,n+1)}\subs U(x,n)$.
Then
there is a first countable regular space $Y$ is $V^Q$
and for each $y\in Y$ there is a neighbourhood base
$\{U'(y,n):n<{\omega}\}$
of $y$ with
$\overline{U'(y,n+1)}\subs U(y',n)$
such that
\begin{enumerate}[(i)]
\item \label{ext} $U'(x,n)\cap X=U(x,n)$ for $x\in X$ and $n<{\omega}$,
\item \label{disj} if $U(x,n)\cap U(y,m)=\empt$
then $U'(x,n)\cap U'(y,m)=\empt$,
\item \label{subs} if $U(x,n)\subs U(y,m)$ then $U'(x,n)\subs U'(y,m)$,
\item  every $A\in \br X;{\omega};\cap V$ has an accumulation point in $Y$.
\end{enumerate}
\end{lemma}

\begin{remark}
If $X$ is 0-dimensional then so is $Y$, and in this case the proof
can be simplified a bit.
\end{remark}

\begin{proof}\prlabel{lm:main}

  \begin{qtheorem}[Hechler, \cite{NyPS}]
   If $W$ is a Cohen generic extension of $V$ then in $W$
there is an almost disjoint family $\bcal\subs \br {\omega};{\omega};$
which refines $\br {\omega};{\omega};\cap V$.
  \end{qtheorem}

Let $\scal\subs \br X;{\omega};\cap V$ be a maximal almost disjoint family in
$V$. By the theorem above for each $S\in \scal$ there is a maximal almost
disjoint family $\tcal_S\subs \br S;{\omega};$ in $V^{\ccal}$
which refines $\br S;{\omega};\cap V$. Then
$\tcal=\cup\{\tcal_S:S\in \scal\}$ is almost disjoint and refines
$\br X;{\omega};\cap V$.

Let
  \begin{displaymath}
\acal=\{A\in \br X;{\omega};: \text{$A$ is closed discrete in $X$}\}.
    \end{displaymath}
Put
\begin{displaymath}
\bcal=\{B\in\tcal:\text{$B$ is closed discrete in $X$}\}.
\end{displaymath}
Then $\bcal$ refines $\acal$.

Then, in $V^{\ccal*R_{\fcal}}$, for each $B\in \bcal$ there is
$D_B\in \br B;{\omega};$ such that for each $x\in X$ and for each $n<{\omega}$
either $D_B\subs^* U(x,n)$ or $D_B\cap U(x,n)$ is finite.

Let $\dcal=\{D_B:B\in \bcal\}$.

Let  $Y=X\cup \{y_D:D\in D\}$, where $y_D$ are new points.

We will define the topology of $Y$ in $V^{\ccal*R_\fcal*\dcal}$
as follows.

Let $d$ be the dominating real introduced by $\dcal$.
For $x\in X$ and $n<{\omega}$ let
\begin{displaymath}
U'(x,n)=U(x,n)\cup\{y_D:D\subs^* U(x,n)\}.
\end{displaymath}

\newcommand{\dvec}{\vec{D}}
\newcommand{\evec}{\vec{E}}

For $D\in \dcal$ let $\dvec$ be a 1-1 enumeration of $D$, and
for $n\in {\omega}$
let
\begin{displaymath}
U'(y_D,n)=\{y_D\}\cup\bigcup\{U'(\dvec(k), d(k)+n):n\le k<{\omega}\}.
\end{displaymath}

The family $\{U'(y,n):y\in Y, n\in {\omega}\}$
clearly satisfies (\ref{ext}), (\ref{disj}) and (\ref{subs}).

We intend to define the topology on $Y$ as the one induced
 by the neighbourhood base
$\{U'(y,n):y\in Y, n\in {\omega}\}$.
Next  we should prove that $\{U'(y,n):y\in Y, n\in {\omega}\}$ is a
neighbourhood base of a topology.

\begin{claim}
If $t\in U'(v,n)$ then there is $m$ such that
$U'(t,m)\subs U'(v,n)$.
\end{claim}

\begin{proof}[Proof of the claim]\prnolabel\mbox{ }\\
{\em Case 1: $v,t\in X$.}\\ Then there is $m$ such that
$U(t,m)\subs U(v,n)$. Then $U'(t,m)\subs U'(v,n)$.\\
{\em Case 2: $t\in X$ and $v=y_D\in Y\setm X$.}\\
Then $t\in U'(\dvec(k), d(k)+n)$ for some  $k\ge n$ and so
$t\in U(\dvec(k), d(k)+n)$. Thus there is $m$
such that $U(t,m)\subs U(\dvec(k),d(k)+n)$ and so
$U'(t,m)\subs U'(\dvec(k),d(k)+n)\subs U'(v,n)$.\\
{\em Case 3: $t=y_D\in Y\setm X$ and $v \in  X$.}\\
Then $D\subs^* U(v,n)$.  Fix $k$ such that
$\dvec(k')\in U(v,n)$ for $k'\ge k$. Pick a function
$g:{\omega}\setm k\to {\omega}$ in $V^{\ccal*R_\fcal}$
such that $U(\dvec(k'),g(k'))\subs U(v,n)$.
Then there is $m\ge k$ such that
$d(m')\ge g(m')$ for $m'\ge m$. Then
$U(\dvec(m'), d(m')+m)\subs U(v,n)$ for $m'\ge m$, hence
$U'(\dvec(m'), d(m')+m)\subs U'(v,n)$ for $m'\ge m$,
 and so
$U'(y_D,m)\subs U'(v,n)$.\\
{\em Case 4:
$t=y_D$ and $v=y_E$ for some $D,E\in \dcal$.}\\
Then $t\in U'(\evec(n'),d(n')+n)$ for some
$n'\ge n$. Since   $\evec(n')\in X$  we can apply Case 3 to get
an  $m$ such that $U'(t,m)\subs U'(\evec(n'),d(n')+n)$
and so $U'(t,m)\subs U'(v,n)$.
\end{proof}

Hence the family
$\{U'(y,n):y\in Y, n\in {\omega}\}$ can be considered as the
neighbourhood base of a topology  on $Y$.

\begin{claim}\label{t2reg}
If $t\notin U'(v,n)$ then there is $m$ such that
$U'(t,m)\cap U'(v,n+1)=\empt$.
\end{claim}

\begin{proof}[Proof of the claim]\prnolabel\mbox{ }\\
{\em Case 1: $v,t\in X$.}\\
Since $\overline{U(v,n+1)}\subs U(v,n)$  there is $m$ such that
$U(t,m)\cap U(v,n+1)=\empt$. Then
$U'(t,m)\cap U'(v,n+1)=\empt$.\\
{\em Case 2: $t=y_D\in Y\setm X$ and $v\in X$.  }\\
Since $y_D\notin U'(v,n)$ we have that $D\cap U'(v,n)$ is finite and so
there is $\ell$ such that
$\{\dvec(i):i\ge \ell\}\cap U(v,n)=\empt$.
Then $\{\dvec(i):i\ge \ell\}\cap\overline{U(v,n+1)}=\empt$ and so
we can find a function $g:{\omega}\setm \ell\to {\omega}$
such that $ U(\dvec(i),g(i))\cap U(v,n+1)=\empt$ for $i\ge \ell$.
Then there is $m\ge \ell$ such that
$d(m')\ge g(m')$ for $m'\ge m$. Then
$U(\dvec(i),d(i))\cap U(v,n+1)=\empt$ and so
$U'(\dvec(i),d(i))\cap U'(v,n+1)=\empt$ as well for $i\ge m$.
Thus $U'(t,m)\cap U'(v,n+1) =\empt$.\\
{\em Case 3: $t\in X$ and $v=y_D\in Y\setm X$.  }\\
Since $D$ is closed discrete and $t\notin U(v,n)$ there is $\ell$ such that
$\overline{U(t,\ell)}\cap \{\dvec(i):i\ge n\}=\empt$.
Fix a function $g:{\omega}\setm n\to {\omega}$ such that
$U(t,\ell)\cap U(\dvec(i),g(i))=\empt$ for $i\ge n$.
There is $k\ge n$ such that $d(i)\ge g(i)$ for $i\ge k$.

Since $t\notin U'(v,n)$ we have $t\notin U(\dvec(i),d(i)+n)$
for $i\ge n$. Thus $t\notin \overline{U(\dvec(i),d(i)+n+1)}$
for $i\ge n$.
Fix $m\ge \ell$ such that
$U(t,m)\cap \cup\{U(\dvec(i),d(i)+n+1):n\le i<\ell\}=\empt$.

Then $U(t,m)\cap U'(v,n+1)=\empt$ and so $U'(t,m)\cap U'(y,n+1)=\empt$
as well.\\
{\em Case 4:
$t=y_D$ and $v=y_E$ for some $D,E\in \dcal$.}\\
Since $D$ and $E$ are closed discrete and $E\cap D$ is finite
there is $\ell<{\omega}$ and a function
$g:{\omega}\setm \ell\to {\omega}$ such that
\begin{displaymath}
\bigcup\{U(\dvec(i),g(i):i\ge \ell\}\cap
\bigcup\{U(\evec(i),g(i)):i\ge \ell\}=\empt.
\end{displaymath}
There is $k\ge \ell$ such that $g(i)\le d(i)$ for $i\ge k$. Then
\begin{displaymath}\tag{$*$}
\bigcup\{U(\dvec(i),d(i):i\ge k\}\cap
\bigcup\{U(\evec(i),d(i)):i\ge k\}=\empt.
\end{displaymath}

Since $t\notin U(\evec(i),d(i)+n)$ for $i\ge n$,  by case 2
for each $i\ge n$ there is $j_i$ such that
$U'(t,j_i)\cap U(\evec(i),d(i)+n+1)=\empt$.
Let $m_0=\max\{{j_i:n\le i<k}\}$.
Then
\begin{displaymath}\tag{$**$}
U'(t,m_0)\cap \{U(\evec(i),d(i)+n+1):n\le i<k\}=\empt.
\end{displaymath}
Let $m=\max\{{m_0,k}\}$.
Then
\begin{displaymath}
U'(t,m)\setm \{t\}\subs \bigcup\{U'(\dvec(i),d(i):i\ge k\}
\end{displaymath}
and so by $(*)$ we have
\begin{displaymath}\tag{$*{*}*$}
U'(t,m)\cap \bigcup\{U(\evec(i),d(i)):i\ge k\}=\empt.
\end{displaymath}
Since
\begin{displaymath}
U'(v,n+1)\setm \bigcup\{U(\evec(i),d(i)):i\ge k\}\subs
\{U(\evec(i),d(i)+n+1):n\le i<k\}
\end{displaymath}
$(**)$ and $(*{*}*)$ together yields
$U'(t,m_0)\cap U'(v,n+1)=\empt$.
\end{proof}

Since $\bigcap \{U'(x,n):n<{\omega}\}=\{x\}$ claim \ref{t2reg} implies
that $X$ is a regular space.
\end{proof}

After proving the lemma we can easily get the theorem.
The poset $P$ is obtained by a finite support iteration
$\<P_{\alpha}:{\alpha}\le {\omega}_1\>$
of length $\oo$,  $P_{{\alpha}+1}=P_{\alpha}*\dot Q_{\alpha}$,
where $Q_{\alpha}$ is the poset $Q$ from lemma \ref{lm:main}
in the model $V^{P_{\alpha}}$.

Let $X$ be a regular first countable space from the ground model.
For each $x\in X$ let $\{U(x,n):n<{\omega}\}$
be a neighbourhood base of $x$ such that
$\overline{U(x,n+1)}\subs U(x,n)$.

We will construct first countable spaces  $X_{\alpha}$
with bases $\{U_{\alpha}(x,n):x\in X_{\alpha},n<{\omega}\}$
satisfying $\overline{U_{{\alpha}}(x,n+1)}\subs U_{\alpha}(x,n)$
such that
\begin{enumerate}[(1)]
\item $X_0=X$, $U_0(x,n)=U(x,n)$
\item $X_{\alpha}, \{U_{\alpha}(x,n):x\in X_{\alpha},
n\in {\omega}\}\in V^{P_{\alpha}}$,
\item \label{ext} $U_{\beta}(x,n)\cap X_{\alpha}=U_{\alpha}(x,n)$ for
${\alpha}<{\beta}$,
$x\in X_{\alpha}$ and $n<{\omega}$,
\item \label{disj} if $U_{\alpha}(x,n)\cap U_{\alpha}(y,m)=\empt$
then $U_{\beta}(x,n)\cap U_{\beta}(y,m)=\empt$ for ${\alpha}<{\beta}$,
\item \label{subs} if $U_{\alpha}(x,n)\subs U_{\alpha}(y,m)$
then $U_{\beta}(x,n)\subs U_{\beta}(y,m)$,
\item  every $A\in \br X_{\alpha};{\omega};\cap V^{P_{\alpha}}$
has an accumulation point in $X_{{\alpha}+1}$.
\end{enumerate}

The construction is straightforward: apply lemma \ref{lm:main}
in successor steps and take
$X_{\alpha}=\bigcup \{X_{\zeta}:{\zeta}<{\alpha}\}$, and
$U_{\alpha}(x,n)=\bigcup\{U_{\zeta}(x,n):x\in X_{\zeta}\}$
for limit ${\alpha}$. Clearly $X_{\alpha}$ and 
$\{U_{\alpha}(x,n):x\in X_{\alpha},n<{\omega}\}$
satisfy the inductive requirements (1)--(6).

Since every countable subset of $Y=X_{\oo}$ appears in some intermediate step
the space $Y$ will be countably compact.

\end{proof}

\end{document}